\documentclass[a4paper,10pt]{article}
\usepackage{hyperref}
\usepackage{cleveref}
\hypersetup{hypertexnames = false, bookmarksdepth = 2, bookmarksopen = true, colorlinks, linkcolor = black, citecolor = black, urlcolor = black, pdfstartview={XYZ null null 1}}

\usepackage{afterpage}
\usepackage{amsfonts}
\usepackage[fleqn, leqno]{amsmath}
\usepackage{amssymb}
\usepackage{amsthm}
\usepackage{booktabs}
\usepackage{dcolumn}
\usepackage{diagbox}
\usepackage{enumitem}
\usepackage{mathtools}
\usepackage{parskip}
\usepackage{pdflscape}
\usepackage{thmtools}
\usepackage{tikz-cd}
\usepackage[colorinlistoftodos, textsize = footnotesize]{todonotes}
\usepackage{xparse}
\usepackage{xspace}

\usepackage[T1]{fontenc}
\usepackage{eucal}
\usepackage{microtype}
\usepackage{gitinfo2}

\newcounter{todocounter}
\DeclareDocumentCommand\addreference{g}{\stepcounter{todocounter}\todo[color = blue!30]{\thetodocounter. Add reference\IfNoValueF{#1}{: #1}}\xspace}
\DeclareDocumentCommand\checkthis{g}{\stepcounter{todocounter}\todo[color = red!50]{\thetodocounter. Check this\IfNoValueF{#1}{: #1}}\xspace}
\DeclareDocumentCommand\fixthis{g}{\stepcounter{todocounter}\todo[color = orange!50]{\thetodocounter. Fix this\IfNoValueF{#1}{: #1}}\xspace}
\DeclareDocumentCommand\expand{g}{\stepcounter{todocounter}\todo[color = green!50]{\thetodocounter. Expand\IfNoValueF{#1}{: #1}}\xspace}

\declaretheoremstyle[
  spaceabove = 3pt,
  spacebelow = 3pt,
]{lecture}
\theoremstyle{lecture}
\newtheorem{theorem}{Theorem}
\newtheorem{corollary}[theorem]{Corollary}
\newtheorem{definition}[theorem]{Definition}
\newtheorem{example}[theorem]{Example}

\newtheorem{lemma}[theorem]{Lemma}
\newtheorem{proposition}[theorem]{Proposition}
\newtheorem{remark}[theorem]{Remark}

\makeatletter
\def\gitfootnote{\gdef\@thefnmark{}\@footnotetext}
\makeatother

\hypersetup{pagebackref}

\mathchardef\mhyphen="2D
\newcommand\dash{\nobreakdash-\hspace{0pt}}

\newcommand\plane[1]{P\textsubscript{#1}}
\newcommand\quadric[1]{Q\textsubscript{#1}}

\newcommand\ab{\ensuremath{\mathrm{ab}}}
\newcommand\bounded{\ensuremath{\mathrm{b}}}
\newcommand\Cyc{\ensuremath{\mathrm{Cyc}}}

\newcommand\Dih{\ensuremath{\mathrm{Dih}}}
\newcommand\dg{\ensuremath{\mathrm{dg}}}
\newcommand\env{\ensuremath{\mathrm{e}}}
\newcommand\Ga{\ensuremath{\mathbf{G}_{\mathrm{a}}}}

\newcommand\Gm{\ensuremath{\mathbf{G}_{\mathrm{m}}}}

\newcommand\PGL{\ensuremath{\mathrm{PGL}}}

\newcommand\SL{\ensuremath{\mathrm{SL}}}
\newcommand\Sym{\ensuremath{\mathrm{Sym}}}

\DeclareMathOperator\Aut{Aut}

\DeclareMathOperator\catmod{mod}
\DeclareMathOperator\CCCC{CC}
\DeclareMathOperator\characteristic{char}
\DeclareMathOperator\coh{coh}

\DeclareMathOperator\derived{\mathbf{D}}
\DeclareMathOperator\Ed{Ed}

\DeclareMathOperator\Ext{Ext}
\DeclareMathOperator\gldim{gldim}
\DeclareMathOperator\Gr{Gr}
\DeclareMathOperator\gr{gr}
\DeclareMathOperator\HH{H}
\DeclareMathOperator\HHHH{HH}

\DeclareMathOperator\Hom{Hom}
\DeclareMathOperator\Ind{Ind}
\DeclareMathOperator\Inj{Inj}

\DeclareMathOperator\Kzero{K_0}
\DeclareMathOperator\Lie{Lie}

\DeclareMathOperator\Ord{Ord}
\DeclareMathOperator\Out{Out}
\DeclareMathOperator\pdim{pdim}
\DeclareMathOperator\Pic{Pic}
\DeclareMathOperator\projdim{projdim}
\DeclareMathOperator\Qd{Qd}
\DeclareMathOperator\qgr{qgr}
\DeclareMathOperator\rk{rk}

\DeclareMathOperator\Tot{Tot}
\DeclareMathOperator\trace{tr}

\title{Hochschild cohomology of noncommutative planes and quadrics}
\author{Pieter Belmans}

\begin{document}
\maketitle


\begin{abstract}
  We give a description of the Hochschild cohomology for noncommutative planes (resp.\ quadrics) using the automorphism groups of the elliptic triples (resp.\ quadruples) that classify the Artin--Schelter regular~$\mathbb{Z}$\dash algebras used to define noncommutative planes and quadrics. For elliptic triples the description of these automorphism groups is due to Bondal--Polishchuk, for elliptic quadruples it is new.
\end{abstract}

\footnotetext{14A22, 16E40, 18E30}

\setcounter{tocdepth}{2}

\section{Introduction}
\label{section:introduction}
Noncommutative algebraic geometry in the sense of Artin--Zhang concerns the abelian categories associated to (not necessarily commutative) graded algebras, giving rise to the category of coherent sheaves on (nonexistent) noncommutative projective varieties. An important class of such algebras is given by the analogues of the (commutative) polynomial ring, and their properties have been axiomatised by Artin--Schelter \cite{MR917738}. A full classification of~3\dash dimensional Artin--Schelter regular algebras has been obtained by Artin--Tate--Van den Bergh \cite{MR1086882,MR1128218}, and it can be seen as the classification of \emph{noncommutative planes} and \emph{noncommutative quadrics}. Indeed, there are precisely two types of graded algebras in this case, distinguished by their Hilbert series, and which correspond to the noncommutative analogues of~$\mathbb{P}^2$ or~$\mathbb{P}^1\times\mathbb{P}^1$.

More recently, Lowen--Van den Bergh developed a deformation theory for abelian categories, generalising that of associative algebras \cite{MR2183254,MR2238922,MR3050709}. They introduced a notion of Hochschild cohomology for abelian categories, which in turn describes the infinitesimal deformations as an abelian category. Based on the computation of the Hochschild cohomology of~$\coh\mathbb{P}^1\times\mathbb{P}^1$ (see \cref{example:hkr-quadric}) it becomes clear that the classification of noncommutative quadrics is incomplete if one only uses graded algebras. There should be 3~degrees of freedom, whereas graded algebras only provide~2~degrees \cite{MR2836401}. In op.~cit.~Van den Bergh generalised the construction of noncommutative quadrics to use~$\mathbb{Z}$\dash algebras, and provides a classification of these in terms of geometric and linear data. For noncommutative planes (where graded algebras do suffice, as in Artin--Tate--Van den Bergh \cite{MR1086882,MR1128218}) the classification using~$\mathbb{Z}$\dash algebras \cite{MR1230966} gives a more streamlined approach which will be used here.

In this paper we describe the Hochschild cohomology of all noncommutative planes and quadrics using geometric techniques. For the commutative case this computation is an easy application of the Hochschild--Kostant--Rosenberg theorem as in \cref{example:hkr-plane,example:hkr-quadric}, but in the noncommutative case no such techniques are available. The main observation from \cref{table:hochschild-quadratic,table:hochschild-cubic} is that noncommutative planes and quadrics have \emph{less (infinitesimal) symmetries} or equivalently that they are \emph{more (infinitesimally) rigid}, which explains the \emph{dimension drop} we see happening in \cref{table:hochschild-quadratic,table:hochschild-cubic}, and we give an explicit description of how this dimension drop behaves for all cases in the classification.

\begin{table}[p]
  \centering
  \small
  \begin{tabular}{rp{4.5cm}cc}
    \toprule
       & divisor                            & $\dim_k\HHHH_\ab^1(\qgr A)$ & $\dim_k\HHHH_\ab^2(\qgr A)$ \\
    \midrule
    1. & elliptic curve                     & 0                           & 2 \\
    2. & cuspidal curve                     & 0                           & 2 \\
    3. & nodal curve                        & 0                           & 2 \\
    4. & three lines in general position    & 2                           & 4 \\
    5. & three lines through a point        & 2                           & 4 \\
    6. & conic and line in general position & 1                           & 3 \\
    7. & conic and tangent line             & 1                           & 3 \\
    8. & triple line                        & 5                           & 7 \\
    9. & double line and line               & 3                           & 5 \\
    \midrule
       & commutative plane                  & 8                           & 10 \\
    \bottomrule
  \end{tabular}
  \caption{Hochschild cohomology of noncommutative planes}
  \label{table:hochschild-quadratic}
\end{table}

\begin{table}[p]
  \centering
  \small
  \begin{tabular}{rp{5cm}ccccc}
    \toprule
        & divisor                                       & $\dim_k\HHHH_\ab^1(\qgr A)$ & $\dim_k\HHHH_\ab^2(\qgr A)$ \\
    \midrule
    1.  & elliptic curve                                & 0                           & 3 \\
    2.  & cuspidal curve                                & 0                           & 3 \\
    3.  & nodal curve                                   & 0                           & 3 \\

    4.  & two conics in general position                & 1                           & 4 \\
    5.  & two tangent conics                            & 1                           & 4 \\
    6.  & conic and two lines in a triangle             & 1                           & 4 \\
    7.  & conic and two lines through a point           & 2                           & 5 \\
    8.  & quadrangle                                    & 2                           & 5 \\
    9.  & twisted cubic and a bisecant                  & 1                           & 4 \\
    10. & twisted cubic and a tangent line              & 1                           & 4 \\

    11. & double conic                                  & 3                           & 6 \\
    12. & two double lines                              & 3                           & 6 \\
    13. & double line and two lines in general position & 3                           & 6 \\
    \midrule
        & commutative quadric                           & 6                           & 9 \\
    \bottomrule
  \end{tabular}
  \caption{Hochschild cohomology of noncommutative quadrics}
  \label{table:hochschild-cubic}
\end{table}

In \cref{subsection:hh-abelian,subsection:hh-varieties-algebras} we recall the various notions of Hochschild cohomology of algebras, varieties and abelian categories and their properties required for this paper. The translation between smooth projective varieties and associative algebras will use exceptional objects, whose properties we quickly recall in \cref{subsection:exceptional}. In \cref{subsection:as-regular} we describe the technicalities regarding~$\mathbb{Z}$\dash algebras which are necessary to define noncommutative planes and quadrics and describe their classification. In \cref{subsection:segre} we give an overview of a completely classical but not so well-known classification of base loci of pencils of quadrics, which arise as the point schemes of noncommutative quadrics.

The method of computing the Hochschild cohomology is given in \cref{section:hh}: it is based on the Lefschetz trace formula for Hochschild cohomology as recalled in \cref{subsection:euler-characteristic}, and the description of the first Hochschild cohomology of a finite-dimensional algebra as the Lie algebra of its outer automorphisms in \cref{subsection:Lie-Out}. In \cref{subsection:identifying-aut} it is explained how this group of outer automorphisms of the finite-dimensional algebra can be related to the automorphism groups of the corresponding~$\mathbb{Z}$\dash algebra, and hence by the geometric classification of noncommutative planes (resp.~quadrics) by Bondal--Polishchuk (resp.~Van den Bergh) we can reduce the problem to a purely geometric question regarding automorphism groups of certain (possibly singular, possibly reducible, possibly nonreduced) curves of arithmetic genus~1.

In \cref{section:aut} we then give a description of the automorphism groups necessary to apply \cref{corollary:out}. For noncommutative planes this description is given in \cite[table~1]{MR1230966}, and recalled in \cref{table:hochschild-quadratic-automorphisms}. For noncommutative quadrics this description is new, and the details for the computations are given in \cref{subsection:noncommutative-quadrics}.

Another class of finite-dimensional algebras arising from algebraic geometry for which the Hochschild cohomology is completely known is that of weighted projective lines \cite{MR1659132}. Remark that in that case it turns out that there is \emph{no} dependence on the parameters (i.e.~the location of the stacky points).

\paragraph{Conventions} Throughout we will assume~$k$ to be an algebraically closed field of characteristic not equal to~$2,3$. When comparing Hochschild cohomology to Poisson cohomology work over the complex numbers.

\paragraph{Acknowledgements}
The author would like to thank Shinnosuke Okawa, Brent Pym, Fumiya Tajima and Michel Van den Bergh for interesting discussions.

The author was supported by a Ph.D.~fellowship of the Research Foundation--Flanders (FWO).

\section{Preliminaries}
\label{section:preliminaries}
In this section we recall the definition of Hochschild cohomology for different types of objects, and its properties that we will need later.

\subsection{Hochschild cohomology of abelian and differential graded categories}
\label{subsection:hh-abelian}
It is well-known that Hochschild cohomology for associative algebras governs their deformation theory. Its definition has been generalised to arbitrary dg~categories \cite{keller-dih}, which we will quickly recall now. In \cref{subsection:hh-varieties-algebras} we will specialise this definition to associative algebras and smooth projective varieties, and discuss some special properties.

\begin{definition}
  Let~$\mathcal{C}$ be a~$k$\dash linear dg~category. Its \emph{Hochschild complex}~$\CCCC_\dg^\bullet(\mathcal{C})$ is the~$\Tot^\prod$ of the complex whose term in degree~$p\geq 0$ is given by
  \begin{equation}
    \prod_{C_0,\dotsc,C_p}\Hom_k\left(\mathcal{C}(C_{p-1},C_p)^\bullet\otimes_k\mathcal{C}(C_{p-2},C_{p-1})^\bullet\otimes_k\dotsc\otimes_k\mathcal{C}(C_0,C_1)^\bullet,\mathcal{C}(C_0,C_p)^\bullet \right)
  \end{equation}
  with differentials as for the Hochschild complex of an associative algebra, and which is zero in degree~$p\leq -1$. The cohomology of this complex is the \emph{Hochschild cohomology of~$\mathcal{C}$}, and will be denoted~$\HHHH_\dg^\bullet(\mathcal{C})$.
\end{definition}

We will also need the Hochschild cohomology of an abelian category~$\mathcal{A}$. This is \emph{not} the Hochschild cohomology of this abelian category considered as a~$k$\dash linear category giving rise to a dg~category in degree~0. Rather, the definition is as follows.

\begin{definition}
  Let~$\mathcal{A}$ be a~$k$\dash linear abelian category. Its \emph{Hochschild cohomology}~$\HHHH_\ab^\bullet(\mathcal{A})$ is defined as the Hochschild cohomology of the dg~category associated to the~$k$\dash linear subcategory of injective objects of the Ind-completion of~$\mathcal{A}$, i.e.
  \begin{equation}
    \HHHH_\ab^\bullet(\mathcal{A})\coloneqq\HHHH_\dg^\bullet(\Inj\Ind\mathcal{A}).
  \end{equation}
\end{definition}

This definition makes sense because for any abelian category we have that~$\Ind\mathcal{A}$ is a Grothendieck abelian category, hence has enough injective objects.

In \cite[theorem~6.1]{MR2183254} it is shown that the Hochschild cohomology of abelian categories agrees with the Hochschild cohomology (as a dg category) of a dg~enhancement of the (bounded) derived category of an abelian category in the following way.
\begin{proposition}[Lowen--Van den Bergh]
  \label{proposition:lowen-vandenbergh}
  Let~$\mathcal{A}$ be an abelian category. Let~$\derived_\dg^\bounded(\mathcal{A})$ be the dg enhancement of~$\derived^\bounded(\mathcal{A})$ given by the full subcategory (of the dg category of complexes in~$\Ind\mathcal{A}$ with injective components) on the left bounded complexes of injectives with bounded cohomology inside~$\mathcal{A}$. Then there exists an isomorphism
  \begin{equation}
    \HHHH_\ab^\bullet(\mathcal{A})\cong\HHHH_\dg^\bullet(\derived_\dg^\bounded(\mathcal{A})).
  \end{equation}
\end{proposition}
This allows us to prove properties of the Hochschild cohomology using purely algebraic techniques, in particular \cref{proposition:HH1-Lie-Out}, by taking the derived equivalence of our category of interest~$\qgr A$ with~$\catmod kQ/I$ given by \cref{proposition:full-exceptional-collection}.

\subsection{Hochschild cohomology of smooth projective varieties and finite-dimensional algebras}
\label{subsection:hh-varieties-algebras}
The definition of Hochschild cohomology for a differential graded category as given in \cref{subsection:hh-abelian} agrees with the classical definition if one considers an associative algebra as a dg~category with a single object, and morphisms concentrated in a single degree, Moreover the Hochschild cohomology of a finite-dimensional algebra agrees with the Hochschild cohomology of its derived category.

On the other hand, it is possible to give a geometric definition of the Hochschild cohomology of a sufficiently nice scheme, in various ways \cite{MR704600,MR1390671,MR1940241}. It can be shown that these in the setting of smooth and projective varieties over a field they all agree with the Hochschild cohomology of a dg~enhancement of the derived category of the scheme. In particular we will denote by~$\HHHH^\bullet(X)$ the Hochschild cohomology of~$X$ for any of these definitions.

The following theorem comes in many flavours \cite{MR0142598,MR917209,MR1390671}, we will need the following version \cite[corollary~0.6]{MR1940241}.
\begin{proposition}[Hochschild--Kostant--Rosenberg]
  Let~$X$ be a smooth projective variety over a field of characteristic zero or~$p>\dim X$. Then there exists an isomorphism
  \begin{equation}
    \tag{HKR}
    \label{equation:hkr-cohomology}
    \HHHH^i(X)\cong\bigoplus_{p+q=i}\HH^p(X,\bigwedge\nolimits^q\mathrm{T}_X).
  \end{equation}
\end{proposition}

Then the following calculation gives the Hochschild cohomology of the projective plane and quadric.
\begin{example}
  \label{example:hkr-plane}
  For the projective plane we have
  \begin{equation}
    \begin{aligned}
      \HHHH^0(\mathbb{P}^2)&\cong\HH^0(\mathbb{P}^2,\mathcal{O}_{\mathbb{P}^2}), \\
      \HHHH^1(\mathbb{P}^2)&\cong\HH^0(\mathbb{P}^2,\mathrm{T}_{\mathbb{P}^2}), \\
      \HHHH^2(\mathbb{P}^2)&\cong\HH^0(\mathbb{P}^2,\bigwedge\nolimits^2\mathrm{T}_{\mathbb{P}^2})
      \cong\HH^0(\mathbb{P}^2,\mathcal{O}_{\mathbb{P}^2}(3)).
    \end{aligned}
  \end{equation}
  In particular, the dimensions are~$1$,~$8$ and~$10$ respectively.

  Another way of computing these dimensions, closer to the approach taken for noncommutative planes, is by using that
  \begin{equation}
    \begin{aligned}
      \HHHH^1(\mathbb{P}^2)
      &\cong\HH^0(\mathbb{P}^2,\mathrm{T}_{\mathbb{P}^2}) \\
      &\cong\Lie\Aut(\mathbb{P}^2) \\
      &\cong\Lie\PGL_3 \\
      &\cong\mathfrak{sl}_3
    \end{aligned}
  \end{equation}
  is~8\dash dimensional, combined with \cref{corollary:dim-HH-surfaces} to determine the size of~$\HHHH^2(\mathbb{P}^2)$.
\end{example}

\begin{example}
  \label{example:hkr-quadric}
  For the quadric surface we have
  \begin{equation}
    \begin{aligned}
      \HHHH^0(\mathbb{P}^1\times\mathbb{P}^1)&\cong\HH^0(\mathbb{P}^1\times\mathbb{P}^1,\mathcal{O}_{\mathbb{P}^1\times\mathbb{P}^1}), \\
      \HHHH^1(\mathbb{P}^1\times\mathbb{P}^1)&\cong\HH^0(\mathbb{P}^1\times\mathbb{P}^1,\mathrm{T}_{\mathbb{P}^1\times\mathbb{P}^1}), \\
      \HHHH^2(\mathbb{P}^1\times\mathbb{P}^1)&\cong\HH^0(\mathbb{P}^1\times\mathbb{P}^1,\bigwedge\nolimits^2\mathrm{T}_{\mathbb{P}^1\times\mathbb{P}^1})
      \cong\HH^0(\mathbb{P}^1\times\mathbb{P}^1,\mathcal{O}_{\mathbb{P}^1\times\mathbb{P}^1}(2,2)).
    \end{aligned}
  \end{equation}
  In particular, the dimensions are~$1$,~$6$ and~$9$ respectively.

  Another way of computing these dimensions, closer to the approach taken for noncommutative quadrics, is by using that
  \begin{equation}
    \begin{aligned}
      \HHHH^1(\mathbb{P}^1\times\mathbb{P}^1)
      &\cong\HH^0(\mathbb{P}^1\times\mathbb{P}^1,\mathrm{T}_{\mathbb{P}^1\times\mathbb{P}^1}) \\
      &\cong\Lie\Aut(\mathbb{P}^1\times\mathbb{P}^1) \\
      &\cong\Lie\left( (\PGL_2\times\PGL_2)\rtimes\Sym_2 \right) \\
      &\cong\mathfrak{sl}_2\oplus\mathfrak{sl}_2
    \end{aligned}
  \end{equation}
  is~6\dash dimensional, combined with \cref{corollary:dim-HH-surfaces}, to determine the dimension of~$\HHHH^2(\mathbb{P}^1\times\mathbb{P}^1)$.
\end{example}

\subsection{Exceptional objects}
\label{subsection:exceptional}
We will quickly recall the notion of an exceptional collection, as we will use this to reduce the description of the Hochschild cohomology of the abelian category~$\qgr A$ to that of a finite-dimensional algebra. For more background one is referred to \cite[\S1]{MR3545926}.

\begin{definition}
  Let~$\mathcal{C}$ be an Hom-finite triangulated category. Then~$C\in\mathcal{C}$ is \emph{exceptional} if
  \begin{equation}
    \Hom_{\mathcal{C}}(C,C[n])\cong
    \begin{cases}
      k & n=0 \\
      0 & n\neq 0.
    \end{cases}
  \end{equation}
  A sequence of exceptional objects~$C_1,\dotsc,C_n$ is an \emph{exceptional collection} if
  \begin{equation}
    \Hom_{\mathcal{C}}(C_i,C_j)=0
  \end{equation}
  whenever~$i>j$. It is \emph{full} if it generates the triangulated category. It is \emph{strong} if
  \begin{equation}
    \Hom_{\mathcal{C}}(C_i,C_j[n])=0
  \end{equation}
  whenever~$i\leq j$ and~$n\neq 0$.
\end{definition}

Then by the following result of Baer and Bondal \cite{MR928291,MR992977} we have can study triangulated categories with full and strong exceptional collections (and more generally tilting objects) using purely algebraic means.
\begin{proposition}
  Let~$\mathcal{C}$ be a triangulated category that admits a full and strong exceptional collection~$C_1,\dotsc,C_n$. Denote~$A$ the endomorphism algebra of~$\bigoplus_{i=1}^nC_i$. Then the functor
  \begin{equation}
    \Hom_{\mathcal{C}}\left( \bigoplus_{i=1}^nC_i,- \right)\colon\mathcal{C}\to\derived^\bounded(A)
  \end{equation}
  is an equivalence of triangulated categories.
\end{proposition}

\begin{remark}
  Hochschild cohomology is not an additive invariant, unlike e.g.~Hochschild homology or algebraic K-theory (for an introduction to additive invariants one is referred to \cite{MR3379910}). Therefore the presence of a full exceptional collection in the derived category is \emph{not} enough to determine the Hochschild cohomology. On the other hand, if~$\mathbb{I}$ denotes an additive invariant, then if~$\mathcal{C}$ has a full exceptional collection of length~$n$ we obtain that
  \begin{equation}
    \mathbb{I}(\mathcal{C})\cong\mathbb{I}(k)^{\oplus n},
  \end{equation}
  i.e.~$\mathbb{I}(\mathcal{C})$ only depends on the length of the exceptional collection and the values that~$\mathbb{I}$ takes for the base field.

  For Hochschild cohomology there will be a subtle dependence on gluing functors, as in \cite{keller-dih,MR3420334}. On the algebraic level this means that the Hochschild cohomology depends heavily on the relations in the quiver, as evidenced by the results in this paper.
\end{remark}

\subsection{Artin--Schelter regular \texorpdfstring{$\mathbb{Z}$}{Z}-algebras}
\label{subsection:as-regular}
In noncommutative algebraic geometry \`a la Artin--Zhang \cite{MR1304753} one studies the analogues of projective varieties by considering noncommutative graded algebras and their ``categories of (quasi)coherent sheaves'' on them, by mimicking the description of~$\coh X$ provided by Serre's theorem. An important class of such noncommutative graded algebras is given by the appropriate analogues of the (commutative) polynomial ring (which defines~$\mathbb{P}^n$), and a suitable definition for these was found by Artin--Schelter \cite{MR917738}.

For the purposes of this paper it will be important to extend the class of algebras from which we will construct noncommutative projective varieties to include~$\mathbb{Z}$\dash algebras \cite[\S4]{MR1230966}. Recall that a~$\mathbb{Z}$\dash algebra is a non-unital associative algebra~$A=\bigoplus_{i,j\in\mathbb{Z}}A_{i,j}$ for which~$A_{i,j}$ is a finite-dimensional vector space, the subalgebra~$A_{i,i}$ is isomorphic to~$k$ for all~$i\in\mathbb{Z}$, and~$A_{i,j}=0$ if~$j>i$. The multiplication law takes the grading into account, in the sense that it is of the form
\begin{equation}
  A_{l,k}\otimes A_{i,j}\to\delta_l^jA_{i,k}.
\end{equation}
The constructions from \cite{MR1304753} can be generalised to this setting \cite[\S2]{MR2836401}, in particular we have an abelian category~$\qgr A$ for every (sufficiently well-behaved)~$\mathbb{Z}$\dash algebra, which is the analogue of~$\coh X$ for a smooth projective variety~$X$. As in op.~cit.~we denote~$P_{i,A}=e_iA$, where~$e_i$ is the local unit associated to~$i\in\mathbb{Z}$, and we denote~$S_{i,A}$ the unique simple quotient of~$P_i$. The categories~$\gr A$ (resp.~$\Gr A$) are the analogues of the categories of graded modules over a graded algebra, and are defined in \cite{MR2836401}

\begin{definition}
  Let~$A$ be a~$\mathbb{Z}$\dash algebra. Then~$A$ is \emph{Artin--Schelter regular} if
  \begin{enumerate}
    \item $A$ is connected, i.e.~$\dim_k A_{i,j}<+\infty$ for all~$i\leq j$,~$A_{i,j}=0$ for all~$i>j$, and~$A_{i,i}\cong k$ for all~$i$;
    \item $\dim_k A_{i,j}$ has growth bounded by a polynomial in~$j-i$;
    \item $\pdim S_{i,A}<+\infty$, and is moreover bounded independently of~$i$;
    \item for all~$i\in\mathbb{Z}$
      \begin{equation}
        \sum_{j,k\in\mathbb{Z}}\dim_k\Ext_{\Gr A}^j(S_{k,A},P_{i,A})=1.
      \end{equation}
  \end{enumerate}
\end{definition}

There is no classification of~3\dash dimensional Artin--Schelter regular~$\mathbb{Z}$\dash algebras as such, but based on the classification of~3\dash dimensional Artin--Schelter regular graded algebras into two classes distinguished by their Hilbert series (or similarly by the shape of the minimal projective resolution of the simple module~$k$) \cite{MR917738} the following definition is taken from \cite[definition~3.1]{MR2836401}.
\begin{definition}
  Let~$A$ be an Artin--Schelter regular~$\mathbb{Z}$\dash algebra. Then
  \begin{enumerate}
    \item $A$ is \emph{3\dash dimensional quadratic} if the minimal resolution of~$S_{i,A}$ is of the form
      \begin{equation}
        0\to P_{i+3,A}\to P_{i+2,A}^{\oplus 3}\to P_{i+1,A}^{\oplus 3}\to P_{i,A}\to S_{i,A}\to 0;
      \end{equation}
    \item $A$ is \emph{3\dash dimensional cubic} if the minimal resolution of~$S_{i,A}$ is of the form
      \begin{equation}
        0\to P_{i+4,A}\to P_{i+3,A}^{\oplus 2}\to P_{i+1,A}^{\oplus 2}\to P_{i,A}\to S_{i,A}\to 0.
      \end{equation}
  \end{enumerate}
\end{definition}
Based on various properties of the graded algebras and their corresponding abelian categories, the category~$\qgr A$ associated to a~3\dash dimensional quadratic (resp.~cubic) algebra~$A$ is called a \emph{noncommutative plane} (resp.~\emph{noncommutative quadric}).

\begin{example}
  \label{example:z-algebra-commutative-plane}
  Let~$B=k[x,y,z]$ be the commutative polynomial ring. One can associate a~$\mathbb{Z}$\dash algebra~$\check{B}$ to it by setting~$\check{B}_{i,j}=B_{j-i}$. Then~$\qgr B\cong\qgr\check{B}\cong\coh\mathbb{P}^2$.
\end{example}

\begin{example}
  \label{example:z-algebra-commutative-quadric}
  The~$\mathbb{Z}$\dash algebra for the category~$\coh\mathbb{P}^1\times\mathbb{P}^1$ is \emph{not} given by mimicking \cref{example:z-algebra-commutative-plane} for the algebra~$k[x,y,z,w]/(xy-zw)$. The trivial module~$k$ has infinite global dimension in this case, so its minimal projective resolution cannot be of the prescribed form. This reasoning actually tells us that the only commutative Artin--Schelter regular graded algebra is the commutative polynomial ring.

  Rather we need to look for a noncommutative graded algebra (or~$\mathbb{Z}$\dash algebra)~$B$ such that~$\qgr B\cong\coh\mathbb{P}^1\times\mathbb{P}^1$. This graded algebra is given by \cite[example~5.1.1]{MR2836401}
  \begin{equation}
    B\coloneqq k\langle x,y\rangle/(x^2y-yx^2,xy^2-y^2x).
  \end{equation}
\end{example}

We will now summarise the main result of \cite{MR1230966,MR2836401}, which is the classification of the Artin--Schelter regular~$\mathbb{Z}$\dash algebras corresponding to noncommutative planes (resp.~quadrics). It gives an equivalence of categories between
\begin{enumerate}
  \item the category of~$\mathbb{Z}$\dash algebras with morphisms being the isomorphisms,
  \item the category of elliptic triples (resp.~quadruples) with morphisms being the isomorphisms (whose precise definition is given in \cref{section:aut}),
  \item the category of linear algebra data that describes the composition law in the quiver associated to the canonical full and strong exceptional collection from \cref{proposition:full-exceptional-collection}.
\end{enumerate}

\begin{remark}
  \label{remark:sierra}
  Observe that isomorphisms of~$\mathbb{Z}$\dash algebras correspond to equivalences of the graded module categories \cite{MR2776790}. We will come back to this in the context of noncommutative planes (resp.~quadrics) in \cref{remark:okawa} (resp.~\cref{remark:translation-principle}).
\end{remark}

In \cref{corollary:out} we will relate the outer automorphisms of the finite-dimensional algebra described in the third part of the comparison to certain automorphisms of the geometric data describing the noncommutative plane (resp.~cubic), which will allow us to give a description of the Hochschild cohomology.

We will change the terminology and notation from \cite[\S5, \S6]{MR1230966} to be consistent with that of \cite{MR2836401}. In that case the linear algebra data describing the~$\mathbb{Z}$\dash algebra of a noncommutative plane is given by the following definition.
\begin{definition}
  A \emph{geometric quadruple} is a quadruple of vector spaces~$(V_0,V_1,V_2,W)$ such that~$\dim_kV_i=3$, and~$W\subseteq V_0\otimes V_1\otimes V_2$ has~$\dim_kW=1$, hence we may assume that~$W=kw$ for some tensor~$w$. We moreover ask that for all~$j=0,1,2$ and for all restrictions~$w|_{v_j}$ for~$v_j\in V_j\setminus\{0\}$ the associated bilinear form~$V_{j-1}\otimes V_{j+1}\to k$ is of rank at least two (where the indices are taken modulo~3).

  We will say that it is \emph{elliptic} if the determinant associated to the restriction~$w_{\phi_j}$ is a cubic in~$\mathbb{P}(V_j)$.
\end{definition}

For the~$\mathbb{Z}$\dash algebra describing a noncommutative quadric the linear algebra data is described as follows.
\begin{definition}
  A \emph{geometric quintuple} is a quintuple of vector spaces~$(V_0,V_1,V_2,V_3,W)$ such that~$\dim_kV_i=2$, and~$W\subseteq V_0\otimes V_1\otimes V_2\otimes V_3$ has~$\dim_kW=1$, hence we may assume that~$W=kw$ for some tensor~$w$. We moreover ask that for all~$j=0,1,2,3$ and for all~$\phi_j\in V_j^\vee\setminus\{0\}$, $\phi_{j+1}\in V_{j+1}^\vee\setminus\{0\}$
  \begin{equation}
    \langle\phi_j\otimes\phi_{j+1},w\rangle\neq 0,
  \end{equation}
  where the indices are taken modulo~4.

  Using \cite[theorem~4.4]{MR2836401} we will say that it is \emph{elliptic} if it is not of the form~$F(A)$ where~$A$ is a linear quadric, in the notation of loc.~cit.

  An isomorphism of quintuples is an isomorphism of the vector spaces~$V_i$ that preserves the tensor~$w$.
\end{definition}

The classification result for noncommutative planes (resp.~quadrics) is now given by a triangle of categorical equivalences. For noncommutative planes it is a combination of \cite[theorem~6.1 and theorem~6.2]{MR1230966}, together with the results of \S4 of~op.~cit.
\begin{proposition}[Bondal--Polishchuk]
  \label{proposition:classification-planes}
  There exist equivalences of categories between
  \begin{enumerate}
    \item elliptic quadratic Artin--Schelter regular~$\mathbb{Z}$\dash algebras;
    \item elliptic triples;
    \item elliptic geometric quadruples;
  \end{enumerate}
  where the morphisms in each category are the isomorphisms of these objects.
\end{proposition}

For noncommutative quadrics it is given in \cite[\S4]{MR2836401}.
\begin{proposition}[Van den Bergh]
  \label{proposition:classification-quadrics}
  There exist equivalences of categories between
  \begin{enumerate}
    \item elliptic cubic Artin--Schelter regular~$\mathbb{Z}$\dash algebras;
    \item elliptic quadruples;
    \item elliptic geometric quintuples;
  \end{enumerate}
  where the morphisms in each category are the isomorphisms of these objects.
\end{proposition}

One way of interpreting the linear algebra data in the classification is by the following description of the derived category of a noncommutative plane (resp.~quadric).
\begin{proposition}
  \label{proposition:full-exceptional-collection}
  Let~$A$ be a quadratic (resp.~cubic) three-dimensional Artin--Schelter regular~$\mathbb{Z}$\dash algebra. Then~$\derived^\bounded(\qgr A)$ admits a full and strong exceptional collection
  \begin{equation}
    \derived^\bounded(\qgr A)=
    \begin{cases}
      \left\langle\, \widetilde{A},\widetilde{A}(1),\widetilde{A}(2)\right\rangle      & \text{$A$ quadratic} \\
      \left\langle\, \widetilde{A},\widetilde{A}(1),\widetilde{A}(2),\widetilde{A}(3)\right\rangle & \text{$A$ cubic}
    \end{cases}
  \end{equation}
  whose structure is described by the quiver
  \begin{equation}
    \begin{tikzpicture}[scale = 2]
      \draw[circle, minimum size = 7pt, inner sep = 0pt]
        (0,0) node (v0) {}
        (1,0) node (v1) {}
        (2,0) node (v2) {};

      \draw (v0) circle (1pt)
            (v1) circle (1pt)
            (v2) circle (1pt);

      \draw[->] (v0) edge [bend left = 45]  node [fill = white] {$x_0$} (v1)
                (v0) edge                   node [fill = white] {$y_0$} (v1)
                (v0) edge [bend left = -45] node [fill = white] {$z_0$} (v1)
                (v1) edge [bend left = 45]  node [fill = white] {$x_1$} (v2)
                (v1) edge                   node [fill = white] {$y_1$} (v2)
                (v1) edge [bend left = -45] node [fill = white] {$z_1$} (v2);
    \end{tikzpicture}
  \end{equation}
  resp.
  \begin{equation}
    \begin{tikzpicture}[scale = 2]
      \draw[circle, minimum size = 7pt, inner sep = 0pt]
        (0,0) node (v0) {}
        (1,0) node (v1) {}
        (2,0) node (v2) {}
        (3,0) node (v3) {};

      \draw (v0) circle (1pt)
            (v1) circle (1pt)
            (v2) circle (1pt)
            (v3) circle (1pt);

      \draw[->] (v0) edge [bend left =  30] node [fill = white]  {$x_0$} (v1)
                (v0) edge [bend left = -30] node [fill = white] {$y_0$} (v1)
                (v1) edge [bend left =  30] node [fill = white] {$x_1$} (v2)
                (v1) edge [bend left = -30] node [fill = white] {$y_1$} (v2)
                (v2) edge [bend left =  30] node [fill = white] {$x_2$} (v3)
                (v2) edge [bend left = -30] node [fill = white] {$y_2$} (v3);
    \end{tikzpicture}
  \end{equation}
  and whose relations in the elliptic case can be obtained from the linear algebra data in \cref{proposition:classification-planes,proposition:classification-quadrics} by considering the tensor~$w$ as a superpotential.
\end{proposition}

In particular we can use the agreement of various notions of Hochschild cohomology as discussed in \cref{subsection:hh-varieties-algebras,subsection:hh-abelian} to get the following description.
\begin{corollary}
  Let~$A$ be a quadratic (resp.~cubic) three-dimensional Artin--Schelter regular~$\mathbb{Z}$\dash algebra. Let~$kQ/I$ be the endomorphism algebra of the full and strong exceptional collection from \cref{proposition:full-exceptional-collection}. Then we have an isomorphism of Gerstenhaber algebras
  \begin{equation}
    \HHHH_\ab^\bullet(\qgr A)\cong\HHHH^\bullet(kQ/I).
  \end{equation}
\end{corollary}

%

\subsection{Segre symbols}
\label{subsection:segre}
The classification of noncommutative planes using geometric data depends on the classification of plane cubic curves \cite[\S6]{MR1230966}. This classification is a classical and well-known result, and will be used without further comment.

The classification of noncommutative quadrics using geometric data depends on the classification of~$(2,2)$\dash divisors on~$\mathbb{P}^1\times\mathbb{P}^1$ \cite[\S4.4]{MR2836401}. This is also classical, but not so well-known. We will give some background to this classification, and in \cref{table:segre-symbols} we will introduce the labeling used in this paper for these divisors.

Such a~$(2,2)$\dash divisor is the base locus of a pencil of quadrics in~$\mathbb{P}^3$, so we wish to classify these base loci. We can represent elements of the pencil as a symmetric~$4\times4$\dash matrix~$\alpha M+\beta N$, where~$M$ is the matrix describing~$\mathbb{P}^1\times\mathbb{P}^1$ in its Veronese embedding (in particular, it is of maximum rank) whilst~$N$ is any other (non-zero) element of the pencil.

\begin{definition}
  Let~$M$ be a square matrix of size~$n+1$. Let~$\lambda_1,\dotsc,\lambda_e$ be the eigenvalues of~$M$, in decreasing order of algebraic multiplicity, and in case of ambiguity increasing order of number of Jordan blocks with said eigenvalue.

  Let~$(m_{j,1},\dotsc,m_{j,f_j})$ be the sizes in decreasing order of the Jordan blocks with eigenvalue~$\lambda_j$ in the Jordan normal form of~$M$. Then the \emph{Segre symbol} of~$M$ is
  \begin{equation}
    \left[ (m_{1,1},\dotsc,m_{1,f_1}),\dotsc,(m_{e,1},\dotsc,m_{e,f_e}) \right].
  \end{equation}
  If~$f_j=1$, then we will write~$m_{j,1}$ instead of~$(m_{j,1})$.
\end{definition}

The following result is the main reason why we are interested in Segre symbols \cite[theorem~8.6.3]{MR2964027}.
\begin{proposition}[Segre]
  The base loci of two pencils of quadrics in~$\mathbb{P}^n$ are projectively equivalent if and only if the Segre symbols coincide and there exists a projective isomorphism of the pencils that preserves the singular quadrics in the pencil.
\end{proposition}

The main example we are interested in is for~$n=3$, in which case there are precisely~14~possible Segre symbols, because there are~14~two\dash dimensional partitions of~4. The case~$[(1,1,1,1)]$ corresponds to the situation where the two quadrics coincide, so we ignore this. Then the possible Segre symbols and a description of the base locus are given in \cref{table:segre-symbols}. They are not ordered based on the degeneration of the coincidence of the eigenvalues, rather they are grouped based on their geometric properties, see also \cref{table:4-types-of-22-divisors}. This mimicks the ordering of the plane cubics used in \cite[\S6]{MR1230966}.

\begin{table}
  \centering
  \begin{tabular}{cp{5cm}c}
    \toprule
                 & divisor                                         & Segre symbol    \\
    \midrule
    \quadric{1}  & elliptic curve                                  & $[1,1,1,1]$     \\
    \quadric{2}  & cuspidal curve                                  & $[3,1]$         \\
    \quadric{3}  & nodal curve                                     & $[2,1,1]$       \\

    \quadric{4}  & two conics in general position                  & $[(1,1),1,1]$   \\
    \quadric{5}  & two tangent conics                              & $[(2,1),1]$     \\
    \quadric{6}  & a conic and two lines in a triangle             & $[2,(1,1)]$     \\
    \quadric{7}  & a conic and two lines intersecting in one point & $[(3,1)]$       \\
    \quadric{8}  & quadrangle                                      & $[(1,1),(1,1)]$ \\
    \quadric{9}  & twisted cubic and a bisecant                    & $[2,2]$         \\
    \quadric{10} & twisted cubic and a tangent line                & $[4]$           \\

    \quadric{11} & double conic                                    & $[(1,1,1),1]$   \\

    \quadric{12} & two double lines                                & $[(2,1,1)]$     \\
    \quadric{13} & a double line and two lines                     & $[(2,2)]$       \\
    \bottomrule
  \end{tabular}
  \caption{Segre symbols for $(2,2)$-divisors on $\mathbb{P}^1\times\mathbb{P}^1$}
  \label{table:segre-symbols}
\end{table}

\section{A description of Hochschild cohomology}
\label{section:hh}
In this section we give a general procedure to determine the dimensions of the Hochschild cohomology groups, the actual computations are performed in \cref{section:aut}. The idea is to use the Euler characteristic of the Hochschild cohomology as in \cref{subsection:euler-characteristic}, which in our case will only involve three terms by \cref{proposition:hh-dim}. As~$\HHHH^0$ will always be one-dimensional in this setting, it suffices to compute~$\HHHH^1$ to determine~$\HHHH^2$, as in \cref{corollary:dim-HH-surfaces}. In \cref{subsection:identifying-aut} we explain how the Lie algebra of outer automorphisms of the finite-dimensional algebra obtained by tilting theory is related to automorphisms of the geometric data used in the classification of noncommutative planes and quadrics.

\subsection{Euler characteristic of Hochschild cohomology}
\label{subsection:euler-characteristic}
A first ingredient in the computation of the Hochschild cohomology of noncommutative planes and quadrics is a Lefschetz type formula describing the Euler characteristic of Hochschild cohomology of a smooth and proper dg category. Recall that any such category admits a Serre functor on its derived category, inducing by functoriality an automorphism~$\mathbb{S}$ of the Hochschild homology.

For any smooth and proper dg category~$\mathcal{C}$ we will denote
\begin{equation}
  \chi(\HHHH^\bullet(\mathcal{C}))=\dim_k\HHHH^{\text{even}}(\mathcal{C})-\dim_k\HHHH^{\text{odd}}(\mathcal{C})=\sum_{i\in\mathbb{Z}}(-1)^i\dim_k\HHHH^i(\mathcal{C})
\end{equation}
the Euler characteristic of Hochschild cohomology, and likewise for the Hochschild homology of~$\mathcal{C}$. We then have the following result \cite[corollary~3.11]{MR3217465}.

\begin{proposition}
  Let~$\mathcal{C}$ be a smooth and proper dg~category. Then
  \begin{equation}
    \chi(\HHHH^\bullet(\mathcal{C}))=\trace\left( \mathbb{S}^{-1}|_{\HHHH_{\text{even}}(\mathcal{C})} \right)-\trace\left( \mathbb{S}^{-1}|_{\HHHH_{\text{odd}}(\mathcal{C})} \right).
  \end{equation}
\end{proposition}

Because~$(-1)^{\dim X}\mathbb{S}^{-1}$ is upper triangular we get the following corollary \cite[example~3.12]{MR3217465}.
\begin{corollary}
  Let~$X$ be a smooth, projective variety. Then
  \begin{equation}
    \label{equation:supertrace-HH}
    \chi(\HHHH^\bullet(X))=(-1)^{\dim X}\chi(\HHHH_\bullet(X)).
  \end{equation}
\end{corollary}
Observe that in general the sum on the left-hand side of \eqref{equation:supertrace-HH} runs from~$0$ to~$2\dim X$, whilst the sum on the right-hand side runs from~$-\dim X$ to~$\dim X$, using the Hochschild--Kostant--Rosenberg theorem.

The main examples we are interested in are~$\mathbb{P}^2$ and~$\mathbb{P}^1\times\mathbb{P}^1$, which admit a full and strong exceptional collection. In general, whenever we have a tilting object (e.g.~the direct sum of a full and strong exceptional collection) we can simplify the right-hand side of \eqref{equation:supertrace-HH}.
\begin{corollary}
  \label{corollary:supertrace-HH-tilting}
  Let~$X$ be a smooth and projective variety with a tilting object. Then
  \begin{equation}
    \label{equation:supertrace-HH-tilting}
    \chi(\HHHH^\bullet(X))=(-1)^{\dim X}\dim_k\HHHH_0(X).
  \end{equation}

  \begin{proof}
    By \cite[theorem~4.1]{MR3034449} we have that~$\HHHH_i(X)$ is concentrated in degree~0.
  \end{proof}
\end{corollary}
Observe that for smooth projective varieties admitting a full and strong exceptional collection we have that~$\dim_k\HHHH_0(X)=\rk\Kzero(X)$ is the number of exceptional objects.

\begin{remark}
  In the situation of \cref{corollary:supertrace-HH-tilting} it is possible to give a purely algebraic proof of \eqref{equation:supertrace-HH-tilting} using \cite[theorem~2.2]{MR1444101}. In this case one uses that the Serre functor acts unipotently on K-theory \cite[lemma~3.1]{MR1230966}, which is an invariant of the derived category hence applies to the finite-dimensional algebra~$kQ/I$.

  It is important to note that the Serre functor acting unipotently \emph{only depends} on the structure of the quiver with relations, hence for any noncommutative plane or quadric we will obtain the same result.
\end{remark}

This already allows for a heuristic interpretation of the moduli spaces in the case of noncommutative surfaces having a full and strong exceptional collection. Observe that we always have~$\HHHH^0(X)\cong k$, so if~$\HHHH^i(X)=0$ for~$i\geq 3$ (i.e.\ that all deformations are unobstructed), we have that
\begin{equation}
  \dim\HHHH^2(X)-\dim\HHHH^1(X)=\rk\Kzero(X)-1
\end{equation}
describes the number of moduli for (noncommutative) deformations of~$X$. This is indeed confirmed in the case of noncommutative planes (resp.~quadrics) by the classification of \cite{MR1086882,MR1230966} (resp.~\cite{MR2836401}). Observe that in both these cases the global dimension of the algebra associated to the full and strong exceptional collection is~2, see \cref{proposition:hh-dim}.

\begin{example}
  The classification of noncommutative projective planes is performed in \cite{MR1230966} in terms of quadratic Artin--Schelter regular~$\mathbb{Z}$\dash algebras. As discussed in \cref{subsection:noncommutative-planes} the geometric data classifying these algebras is an elliptic triple. In the generic case we are considering elliptic curves, and we see that there are~$1+2-1=2$ moduli: one from the~$j$\dash line, two from~$\Pic^3 C$ and we subtract one from the one-dimensional automorphism group of~$C$.
\end{example}

\begin{example}
  The classification of noncommutative quadrics is performed in \cite{MR2836401} in terms of cubic Artin--Schelter regular~$\mathbb{Z}$\dash algebras. As discussed in \cref{subsection:noncommutative-quadrics} the geometric data classifying these algebras is an elliptic quadruple. In the generic case we are considering elliptic curves, and we see that there are~$1+3-1=3$ moduli: one from the~$j$\dash line, three from~$\Pic^2 C$ and we subtract one from the one-dimensional automorphism group of~$C$.
\end{example}
Similar observations apply to the other del Pezzo surfaces by \cite{MR1846352}.

\subsection{First Hochschild cohomology as Lie algebra}
\label{subsection:Lie-Out}
As explained in \cite{MR2043327}, we have that Hochschild cohomology and its structure as a super-Lie algebra is an invariant of the derived equivalence class. In particular, when computing~$\HHHH_\ab^\bullet(\qgr A)$ we can compute this by using tools from the representation theory of finite-dimensional algebras to compute~$\HHHH^\bullet(kQ/I)$.

We are interested in the first Hochschild cohomology group, for which we have the following description \cite{MR2076381,MR2043327}.
\begin{proposition}
  \label{proposition:HH1-Lie-Out}
  Let~$B\coloneqq kQ/I$ be a finite-dimensional algebra. There exists an isomorphism of Lie algebras
  \begin{equation}
    (\HHHH^1(B),[-,-])\cong\Lie\Out(B)=\Lie\Out^0(B).
  \end{equation}
\end{proposition}
Observe that this therefore describes the Lie algebra structure of~$\HHHH_\ab^1(\qgr A)$. It would be interesting to describe both the algebra structure of~$\HHHH_\ab^\bullet(\qgr A)$ and the Lie module structure of~$\HHHH_\ab^{\geq 2}(\qgr A)$. For the algebra structure it is known (using a purely algebraic proof) for the commutative plane and quadric that the cup products of elements in degree~1 generate the degree~2~part. On the other hand, computations suggest that for noncommutative planes and quadrics the algebra structure is usually trivial, i.e.~all cup products are zero, except for a few cases in the classification.


\begin{remark}
  It is important to observe that~$\Out(A)$ is not necessarily an invariant of the derived category, but~$\Out^0(A)$ is. Also, this description is valid independent of the characteristic of~$k$. Because we avoid~$\characteristic k=2,3$ we will have that the~$\Out(A)$ in question are smooth algebraic groups.
\end{remark}

To compute all dimensions of the Hochschild cohomology of a noncommutative plane (resp.~quadric) we will use that the finite-dimensional algebra obtained by tilting is of global dimension~2, which limits the possibly nonzero Hochschild cohomology spaces. Observe that in the commutative case we could have used the Hochschild--Kostant--Rosenberg decomposition from \eqref{equation:hkr-cohomology}.

\begin{proposition}
  \label{proposition:hh-dim}
  Let~$A$ be a quadratic (resp.~cubic) Artin--Schelter regular~$\mathbb{Z}$\dash algebra. Then we have that
  \begin{equation}
    \HHHH_\ab^i(\qgr A)=0
  \end{equation}
  for~$i\geq 3$.

  \begin{proof}
    Let~$B$ be a finite-dimensional algebra. We will denote the \emph{Hochschild cohomology dimension} by
    \begin{equation}
      \HHHH^\bullet\mhyphen\dim B\coloneqq\sup\{n\in\mathbb{N}\mid\HHHH^n(B)\neq 0\}.
    \end{equation}
    Because~$\projdim_{B^\env}B=\gldim B$ we have that
    \begin{equation}
      \label{equation:hhdim-bound}
      \HHHH^\bullet\mhyphen\dim B\leq\gldim B.
    \end{equation}
    We have that~$\derived^\bounded(\qgr A)$ admits the tilting object constructed as the direct sum of the~$A(i)$'s for~$i=0,\ldots,2$ (resp.~$i=0,\ldots,3$).

    In the \emph{quadratic} case it is immediate that~$\gldim B=2$ as~$\gldim B\leq\#Q_0-1$, and~$B$ is not hereditary.

    In the \emph{cubic} case use that there are precisely two relations of length~3 to conclude by a computation that the projective dimension of the simple object associated to~$A$ is~2.
  \end{proof}
\end{proposition}

Hence under the assumption that~$\gldim kQ/I=2$ and that the Serre functor acts unipotently on Hochschild homology, we can compute the dimension as follows:
\begin{corollary}
  \label{corollary:dim-HH-surfaces}
  \begin{equation}
    \dim_k\HHHH_\ab^i(\qgr A)=
    \begin{cases}
      1 & i=0 \\
      \dim_k\Lie\Out^0(kQ/I) & i=1 \\
      \#Q_0+\dim\HHHH^1(kQ/I)-1 & i=2 \\
      0 & i\geq 3
    \end{cases}.
  \end{equation}

  \begin{proof}
    We have that
    \begin{equation}
      \HHHH_\ab^0(\qgr A)\cong\HHHH^0(kQ/I)\cong\mathrm{Z}(kQ/I)\cong k
    \end{equation}
    because the quiver is acyclic and connected.

    The result now follows from \cref{proposition:HH1-Lie-Out}, \cref{proposition:hh-dim} and the proof of \cref{corollary:supertrace-HH-tilting} because the unipotency of the Serre functor only depends on the structure of the Cartan matrix, which is invariant under modifying the relations in the quiver associated to a noncommutative plane (resp.~quadric).
  \end{proof}
\end{corollary}

\begin{remark}
  Observe that one expects a strong correspondence between the Hochschild cohomology of a noncommutative plane (resp.~quadric) and the Poisson cohomology of its semiclassical limit. In particular, working over the complex numbers the dimension formula in \cref{corollary:dim-HH-surfaces} is closely related to the results of \cite{MR2802550}.
\end{remark}


\subsection{Identifying automorphism groups}
\label{subsection:identifying-aut}
The following theorem is the main tool in computing the Hochschild cohomology of noncommutative planes and quadrics, and depends crucially on the classification result for these objects. The lack of a classification is the main obstruction in systematically generalising the description of Hochschild cohomology for more general noncommutative objects, such as the other noncommutative del Pezzo surfaces, or noncommutative~$\mathbb{P}^3$'s.

To a~$\mathbb{Z}$\dash algebra we can associate its \emph{truncation}. This is a finite-dimensional algebra, which is used in the classification of noncommutative planes and quadrics. If~$A$ is an elliptic quadratic (resp.~cubic) Artin--Schelter regular~$\mathbb{Z}$\dash algebra we define the \emph{truncation}~$A_\sigma$ as
\begin{equation}
  A_\sigma\coloneqq\bigoplus_{i,j=0}^nA_{i,j}
\end{equation}
where~$n=3$ (resp.~$n=4$). It has the following easy but important interpretation in terms of the exceptional collection of \cref{proposition:full-exceptional-collection}.
\begin{lemma}
  \label{lemma:truncated-algebra-interpretation}
  There exists an isomorphism
  \begin{equation}
    A_\sigma\cong kQ/I,
  \end{equation}
  where~$kQ/I$ is the algebra associated to the full and strong exceptional collection from \cref{proposition:full-exceptional-collection}.

  \begin{proof}
    Because~$A$ is generated in degree~1 we can associate the arrows in the quiver with a basis for the components~$A_{i,i+1}$, and the obvious morphism induced from this is an isomorphism by a dimension count.
  \end{proof}
\end{lemma}

We will write~$\mathcal{T}$ (resp.~$\mathcal{Q}$) for the category of elliptic triples (resp.~quadruples) as introduced in \cref{subsection:as-regular}. Then for an object~$(C,\mathcal{L}_0,\mathcal{L}_1)$ we will consider the automorphism group~$\Aut_{\mathcal{T}}(C,\mathcal{L}_0,\mathcal{L}_1)$ (and likewise for elliptic quadruples).

\begin{theorem}
  \label{theorem:aut}
  Let~$A$ be an elliptic~3\dash dimensional quadratic (resp.~cubic)~$\mathbb{Z}$\dash algebra. Let~$A_\sigma$ be the truncation of~$A$ as in \cref{lemma:truncated-algebra-interpretation}. Let~$(C,\mathcal{L}_0,\mathcal{L}_1)$ (resp.~$(C,\mathcal{L}_0,\mathcal{L}_1,\mathcal{L}_2)$) be the elliptic triple (resp.~quadruple) associated to~$A$. Then
  \begin{equation}
    \label{equation:aut-P2}
    \Aut_k(A_\sigma)\cong\Aut_{\mathcal{T}}(C,\mathcal{L}_0,\mathcal{L}_1)
  \end{equation}
  resp.
  \begin{equation}
    \label{equation:aut-P1xP1}
    \Aut_k(A_\sigma)\cong\Aut_{\mathcal{Q}}(C,\mathcal{L}_0,\mathcal{L}_1,\mathcal{L}_2).
  \end{equation}

  \begin{proof}
    For noncommutative planes it is shown in \cite[proposition~5.8]{MR1230966} that the category~$\widetilde{\Ord}_3$ is equivalent to the category of non-degenerate quantum determinants~$\widetilde{\Qd}_3$. This first category is precisely the category of path algebras on~3~vertices with relations, with morphisms being the isomorphisms. We can moreover restrict ourselves now to the subcategory~$\widetilde{\Qd}(3,3,3)$ of~$\widetilde{\Qd}_3$, consisting of non-degenerate quantum determinants of the form~$\varphi\colon U\otimes_kV\otimes_kW\to k$ where~$\dim_kU=\dim_kV=\dim_kW=3$.

    In \cite[theorem~6.2]{MR1230966} it is shown that the subcategory~$\Ed$ of~$\widetilde{\Qd}(3,3,3)$ of tensors for which the determinants of the restriction of~$\varphi$ to~$U$, $V$ or~$W$ are cubics in~$\mathbb{P}^2$ (i.e.~we are in the elliptic case) is equivalent to the category of elliptic triples~$\mathcal{T}$. In particular, the automorphism groups for the~$\mathbb{Z}$\dash algebras for noncommutative planes are identified with the automorphism groups of the associated finite-dimensional algebra, as in \eqref{equation:aut-P2}.

    Similarly for noncommutative quadrics, we can use \cite[corollaries~4.3 and~4.5]{MR2836401} for the relation between finite-dimensional algebras and noncommutative quadrics, and \cite[corollary~4.32]{MR2836401} for the correspondence between automorphisms of the finite-dimensional algebras with the automorphisms of the elliptic quadruples.
  \end{proof}
\end{theorem}

\begin{remark}
  The structure of the objects in the subcategory~$\widetilde{\Qd}(3,3,3)$ in the proof of \cref{theorem:aut} rather corresponds to the exceptional collection~$\mathcal{O}_{\mathbb{P}^2},\mathrm{T}_{\mathbb{P}^2}(-1),\mathcal{O}_{\mathbb{P}^2}(1)$ on~$\mathbb{P}^2$. For the actual computation of the Hochschild cohomology this does not matter.
\end{remark}

To apply \cref{proposition:HH1-Lie-Out} we need to understand the outer automorphisms of the endomorphism algebra of the full and strong exceptional collection. This is achieved by the following corollary, where we denote~$\Aut_{\mathcal{L}_0,\mathcal{L}_1}(C)$ for the subgroup of~$\Aut(C)$ preserving~$\mathcal{L}_0$ and~$\mathcal{L}_1$ (and likewise for an elliptic quadruple).
\begin{corollary}
  \label{corollary:out}
  In the notation of \cref{theorem:aut} we have that
  \begin{equation}
    \label{equation:out-P2}
    \Out_k(kQ/I)\cong\Aut_{\mathcal{L}_0,\mathcal{L}_1}(C)
  \end{equation}
  resp.
  \begin{equation}
    \label{equation:out-P1xP1}
    \Out_k(kQ/I)\cong\Aut_{\mathcal{L}_0,\mathcal{L}_1,\mathcal{L}_2}(C).
  \end{equation}

  \begin{proof}
    By \cref{lemma:truncated-algebra-interpretation} we can replace the truncated algebra~$A_\sigma$ by the endomorphism algebra of the full and strong exceptional collection.

    Now we need to determine the role of the inner automorphisms in the description \eqref{equation:aut-P2} (resp.~\eqref{equation:aut-P1xP1}). The inner automorphisms for~$kQ/I$ are given by conjugation with an element of the form~$\sum_{i=1}^n\alpha_ie_i$, with~$\alpha_i\in k^\times$ and~$e_i$ the idempotents corresponding to the vertices, for~$n=3$ (resp.~$n=4$). These give an inner automorphism group of the form~$k^\times\times k^\times$ (resp.~$k^\times\times k^\times\times k^\times$).

    In the description of the automorphism of elliptic triples (resp.~quadruples), these are precisely the automorphisms of the pairs (resp.~triples) of line bundles: we have that~$\Aut(\mathcal{L}_i)=k^\times$. So we obtain a short exact sequence
    \begin{equation}
      0\to k^\times\times k^\times\to\Aut_{\mathcal{T}}(C,\mathcal{L}_0,\mathcal{L}_1)\to\Aut_{\mathcal{L}_0,\mathcal{L}_1}(C)\to 0
    \end{equation}
    resp.
    \begin{equation}
      0\to k^\times\times k^\times\to\Aut_{\mathcal{Q}}(C,\mathcal{L}_0,\mathcal{L}_1,\mathcal{L}_2)\to\Aut_{\mathcal{L}_0,\mathcal{L}_1,\mathcal{L}_2}(C)\to 0
    \end{equation}
    Hence we have the description of the outer automorphisms of~$kQ/I$ as in \eqref{equation:out-P2} (resp.~\eqref{equation:out-P1xP1}).
  \end{proof}
\end{corollary}
In \cref{section:aut} we describe the automorphism groups in \eqref{equation:out-P2} and \eqref{equation:out-P1xP1}. For noncommutative planes this description is already available in~\cite[table~6.1]{MR1230966}, whilst for noncommutative quadrics it is new.

\subsection{Hochschild cohomology of Artin--Schelter regular algebras}
In the graded case one can also compute the Hochschild cohomology of~$A$ \emph{as an algebra} \cite[\S1.5]{MR1217970}. In the quadratic case the Hochschild (co)homology is computed for skew polynomial algebras in \cite{MR1252939}, and for the generic Sklyanin algebra in \cite{MR1291019}. In the generic Sklyanin-like cubic case there is \cite{MR2071658}.

Indeed, using Poincar\'e--Van den Bergh duality we conclude that the Hochschild cohomology is just dual to the Hochschild homology \cite{MR1443171}. One could ask whether there is a connection between the Hochschild cohomology as an algebra, and the Hochschild cohomology of its associated abelian category. The following examples shows that there is no straightforward connection, by exhibiting two quadratic Artin--Schelter regular algebras, for which the zeroth Hochschild cohomology is not the same, but the abelian categories are equivalent.

\begin{example}
  Consider
  \begin{equation}
    \begin{aligned}
      A&\coloneqq k[x,y,z] \\
      B&\coloneqq k_{\mathbf{q}}[x,y,z]
    \end{aligned}
  \end{equation}
  where~$k_{\mathbf{q}}[x,y,z]$ is the skew polynomial ring associated to the coefficient matrix~$\mathbf{q}=(q_{i,j})$ for which~$q_{i,j}=q_{j,i}^{-1}$. If~$q_{0,1}q_{1,2}q_{0,2}^{-1}=1$, then the point variety of~$B$ is isomorphic to~$\mathbb{P}^2$ \cite{MR3527537}, and in particular~$\qgr A\cong\qgr B$.

  But if the parameters~$q_{i,j}$ are not roots of unity, then the ring is not finite over its center \cite{MR1624000}. In particular we have that
  \begin{equation}
    \HHHH^0(A)=\mathrm{Z}(A)=A\not\cong\HHHH^0(B)
  \end{equation}
\end{example}
This subtle dependence on the choice of the presentation disappears when considering the abelian category~$\qgr A$.

Also, the Hochschild homology of the algebra~$B$ lives in degrees~$0,\ldots,3$ and again depends on the choice of parameters, whereas the Hochschild homology of~$\derived^\bounded(\qgr B)$ is computed from the exceptional collection, and hence isomorphic to~3 (resp.~4) copies of~$k$ \cite{MR3379910}. Likewise there is no immediate connection for the usual commutative polynomial ring~$A$, or for other examples from the classification of Artin--Schelter regular algebras.

\section{Automorphisms of 3-dimensional Artin--Schelter regular \texorpdfstring{$\mathbb{Z}$}{Z}-algebras}
\label{section:aut}
We now describe the automorphism groups which are necessary to use \cref{corollary:out}.
\subsection{Noncommutative planes}
\label{subsection:noncommutative-planes}
The classification of three-dimensional quadratic Artin--Schelter regular algebras from \cref{proposition:classification-planes} by the classification of elliptic (and linear) triples. We are mostly interested in the elliptic case here, as the linear case corresponds to the commutative projective plane, for which one can appeal to Hochschild--Kostant--Rosenberg to compute the Hochschild cohomology as in \cref{example:hkr-plane}.

The notion of elliptic triple used in \cref{proposition:classification-planes} is defined as follows.
\begin{definition}
  An \emph{elliptic triple} is a triple~$(C,\mathcal{L}_0,\mathcal{L}_1)$ where
  \begin{enumerate}
    \item $C$ is a curve;
    \item $\mathcal{L}_0$ and~$\mathcal{L}_1$ are non-isomorphic very ample line bundles of degree 3;
  \end{enumerate}
  such that~$\deg(\mathcal{L}_0|_{C_i})=\deg(\mathcal{L}_1|_{C_i})$ for each irreducible component~$C_i$ of~$C$.
\end{definition}

The curve~$C$ can be interpreted as the \emph{point scheme}, which is a fine moduli space for point modules. We will not need this interpretation, but it is interesting to observe that the infinitesimal automorphisms of the noncommutative object are intimately related to the automorphisms of an associated commutative object.

To apply \cref{corollary:out} for the computation of~$\HHHH^\bullet(\qgr A)$ where~$A$ is an elliptic quadratic Artin--Schelter regular~$\mathbb{Z}$\dash algebra describing a noncommutative plane we need to compute the algebraic group~$\Aut_{\mathcal{L}_0,\mathcal{L}_1}(C)$, where~$(C,\mathcal{L}_0,\mathcal{L}_1)$ is the elliptic triple associated to~$A$ by \cref{proposition:classification-planes}. This description can be found in \cite[table~1]{MR1230966} and is recalled in \cref{table:hochschild-quadratic-automorphisms}. Reading off the dimension of the algebraic groups leads to \cref{table:hochschild-quadratic}.

\begin{table}[tp]
  \renewcommand{\arraystretch}{1.4}
  \centering
  \begin{tabular}{rp{5cm}cc}
    \toprule
    & divisor
    & $\Pic^0C$
    & $\Aut_{\mathcal{L}_0,\mathcal{L}_1}(C)$
    \\
    \midrule
    \plane{1}
    & elliptic curve
    & $C$
    & $(\mathbb{Z}/3\mathbb{Z})^{\oplus2}$
    \\
    \plane{2}
    & cuspidal cubic
    & $\Ga$
    & $1$
    \\
    \plane{3}
    & nodal cubic
    & $\Gm$
    & $\mathbb{Z}/3\mathbb{Z}$
    \\
    \plane{4}
    & three lines in general position
    & $\Gm$
    & $\Gm^2\rtimes\Cyc_3$
    \\
    \plane{5}
    & three concurrent lines
    & $\Ga$
    & $\Ga^2\rtimes\Sym_3$
    \\
    \plane{6}
    & conic and a line
    & $\Gm$
    & $\Gm$
    \\
    \plane{7}
    & conic and a tangent line
    & $\Ga$
    & $\Ga$
    \\
    \plane{8}
    & triple line
    & $\Ga$
    & $\Ga^2\rtimes\SL_2$
    \\
    \plane{9}
    & double line and a line
    & $\Ga$
    & $\Ga\times(\Ga\rtimes\Gm)$
    \\
    \bottomrule
  \end{tabular}
  \caption{Automorphism groups of elliptic triples}
  \label{table:hochschild-quadratic-automorphisms}
\end{table}

The main observation is the cohomology jump, arising from the fact that the commutative plane has more symmetries than any of the noncommutative planes. The generic quadratic Artin--Schelter regular algebra is of type \plane{1}, whose associated finite-dimensional algebra has only finitely many outer automorphisms, hence the associated noncommutative plane has no infinitesimal automorphisms. Allowing the elliptic curve to degenerate further and further will usually increase the size of the automorphism group.

\begin{remark}
  A similar drop in the size of the automorphism group can be observed for the automorphisms of ``the affine complement of the point scheme~$C$''. If~$A$ is a Sklyanin algebra associated to a point of infinite order it is claimed in \cite[proposition~2.10]{MR2303228} that the automorphism group of this~$k$\dash algebra is finite. Remark that in \cref{table:hochschild-quadratic-automorphisms} all Sklyanin algebras are of type \plane{1}, and that the order of the translation does not influence the outer automorphism group of the algebra.
\end{remark}

\begin{remark}
  \label{remark:okawa}
  Observe that we are mostly interested in automorphisms of the algebra~$kQ/I$. These should somehow be related to autoequivalences of~$\qgr A$, but there is no result relating the autoequivalences of an abelian category to its first Hochschild cohomology.

  If~$A$ is a quadratic Artin--Schelter regular~$\mathbb{Z}$\dash algebra the correspondence from \cref{remark:sierra} can be extended to also include equivalences of the quotient category \cite[corollary~A.10]{1411.7770v1}, i.e.~equivalences of~$\qgr A$ can be lifted back to equivalences of~$\gr A$. Such a result is not true for cubic Artin--Schelter regular algebras as we will explain in \cref{remark:translation-principle}.
\end{remark}

\subsection{Noncommutative quadrics}
\label{subsection:noncommutative-quadrics}
The classification of three-dimensional cubic Artin--Schelter regular algebras from \cref{proposition:classification-quadrics} is described by the classification of elliptic (and linear) quadruples. We are mostly interested in the elliptic case here, as the linear case corresponds to the commutative quadric surface, for which one can appeal to Hochschild--Kostant--Rosenberg to compute the Hochschild cohomology as in \cref{example:hkr-quadric}.

The notion of elliptic quadruple used in \cref{proposition:classification-planes} is defined as follows.
\begin{definition}
  An \emph{elliptic quadruple} is a quadruple~$(C,\mathcal{L}_0,\mathcal{L}_1,\mathcal{L}_2)$ where
  \begin{enumerate}
    \item $C$ is a curve;
    \item $(\mathcal{L}_0,\mathcal{L}_1)$ and~$(\mathcal{L}_1,\mathcal{L}_2)$ embed~$C$ into~$\mathbb{P}^1\times\mathbb{P}^1$;
  \end{enumerate}
  such that~$\deg(\mathcal{L}_0|_{C_i})=\deg(\mathcal{L}_2|_{C_i})$ for each irreducible component~$C_i$ of~$C$, and moreover~$\mathcal{L}_0\not\cong\mathcal{L}_2$.
\end{definition}

\begin{remark}
  In \cite{MR2836401} these are called admissible quadruples, but to make the terminology consistent with \cite{MR1230966} we have adapted the definition of ellipticity to already include regularity and non-prelinearity.
\end{remark}

Whereas the full case-by-case description of elliptic triples for noncommutative planes is in the literature, the analogous description for elliptic quadruples is not. The curve~$C$ in an elliptic quadruple is a~$(2,2)$\dash divisor on~$\mathbb{P}^1\times\mathbb{P}^1$, for which we can use Segre symbols to distinguish these. Now to perform this case-by-case analysis, we will distinguish the four types divisors, using the notation from \cref{table:segre-symbols}, as given in \cref{table:4-types-of-22-divisors}.

\begin{table}[h!]
  \centering
  \begin{tabular}{r|cc}
                & reduced                                                                                    & nonreduced \\
    \midrule
    irreducible & \quadric{1}, \quadric{2}, \quadric{3}                                                      & \quadric{11} \\
    reducible   & \quadric{4}, \quadric{5}, \quadric{6}, \quadric{7}, \quadric{8}, \quadric{9}, \quadric{10} & \quadric{12}, \quadric{13}
  \end{tabular}
  \caption{4 types of divisors}
  \label{table:4-types-of-22-divisors}
\end{table}

\paragraph{Reduced and irreducible} These cases (corresponding to cases \quadric{1}, \quadric{2} and \quadric{3}) are covered explicitly in \cite{MR1230966} for the case of noncommutative planes, and the proof can be adapted immediately to the case of~$\mathbb{P}^1\times\mathbb{P}^1$, using~2\dash torsion rather than~3\dash torsion in~$\Pic^0C$.

\paragraph{Reduced, but reducible}
Cases \quadric{4} up to \quadric{10} are of this form.

For cases \quadric{4}, \quadric{6}, \quadric{8} and \quadric{9} we have that~$\Pic^0C\cong\Gm$ because the singularities are all nodal. We can interpret~$\Pic^0C$ as the gluing data for a line bundle.

\begin{remark}
  Cases \quadric{4} and \quadric{9} are isomorphic as abstract schemes, but they differ by their embedding into~$\mathbb{P}^3$ and hence by the properties of the line bundles~$\mathcal{L}_i$.
\end{remark}

In each case the normalisation consists of the appropriate number of copies of~$\mathbb{P}^1$, and to compute the automorphism groups one has to fix the preimages of the singularities on each~$\mathbb{P}^1$, so the automorphism group of each component is the subgroup~$\Gm$ of~$\PGL_2$. Then it is possible to either permute the components or permute the preimages of the singularities, which gives a semidirect product of the torus with a finite group.

We will explicitly describe the action of~$\Aut C$ on~$\Pic C$ for case \quadric{4}, the others are similar. The finite group in this case is~$\Dih_2=\mathbb{Z}/2\mathbb{Z}\times\mathbb{Z}/2\mathbb{Z}$. One copy of~$\mathbb{Z}/2\mathbb{Z}$ exchanges the components, the other switches the preimages (and therefore each component of the normalisation is flipped too). So the action on~$\Gm^2$ of the first copy exchanges the factors, whilst the second copy inverts both~$\Gm$'s.

The copies of~$\Gm$ in~$\Aut C$ act by rescaling~$\Pic^0C$ taking into account the partial degrees of the line bundle, i.e.~$\lambda\in\Gm\subseteq\Aut C$ corresponding to component~$i$ will rescale an element of~$\Pic^0C\oplus\mathbb{Z}^{\oplus 4}$ by~$\lambda^{k_i}$, where~$k_i$ is the partial degree on the component~$i$.

Finally, the first copy of~$\mathbb{Z}$ permutes the copies of~$\mathbb{Z}$ in~$\Pic C$ and inverts~$\Pic^0C$, whilst the other only inverts~$\Pic^0C$.

For cases \quadric{5}, \quadric{7} and \quadric{10} have~$\Pic^0C\cong\Ga$, which can now be interpreted as the tangent space.

\begin{remark}
  Cases \quadric{5} and \quadric{10} are isomorphic as abstract schemes, but their embeddings in~$\mathbb{P}^1\times\mathbb{P}^1$ are different.
\end{remark}

The computation for \quadric{5} (and hence \quadric{10}) is analogous to that of~\plane{7}, except that the line bundles~$\mathcal{L}_i$ will be different. To describe~$\Aut C$ we take the normalisation which has~2~$\mathbb{P}^1$'s, where~1 point is fixed on each copy. So the automorphisms of this configuration are~$(\Ga\rtimes\Gm)^2\rtimes\Sym_2$. Because the singularity is \emph{not} an ordinary multiple point we have to take into account that an automorphism needs to act in the same way on the two tangent spaces at the common point.

Written in affine coordinates an element~$(a_i,b_i)\in\Ga\rtimes\Gm$ acts on~$x\in\mathbb{P}^1\setminus\{\infty\}$ as
\begin{equation}
  x\mapsto a_ix+b_i
\end{equation}
where the point~$\infty$ is the point being fixed. Then for the point at~$\infty$ we have that the action is
\begin{equation}
  x^{-1}\mapsto(a_ix+b_i)^{-1}=x^{-1}(a_i+b_ix^{-1})^{-1}=a_ix^{-1}-b_ia^{-1}x^{-2}
\end{equation}
because infinitesimally we have that~$x^{-2}=0$. So the condition is that~$a_1=a_2$.

For case \quadric{5}, as in the case of \plane{7} we get that fixing the first line bundle (whose~$\Pic^0C$ we take to be~$0$) implies taking the diagonal of the two copies of~$\Ga$. By the choice of~$\mathcal{L}_0$ the copy of~$\Gm$ is preserved. Fixing~$\mathcal{L}_1$ removes this copy. The copy of~$\Sym_2$ is preserved throughout.

For case \quadric{10} on the other hand, we have that the degrees of the line bundles are~$(2,0)$ (resp.~$(1,1)$). Taking line bundle of degree~$(1,1)$ the first step is the same as in the previous paragraph, but the copy of~$\Sym_2$ will disappear in the second step.

The remaining case \quadric{7} is similar to \plane{5}. Again the singularity is not an ordinary multiple point, so to describe the automorphism group we need to take the tangent space into account. If we let~$v_i$ be the vectors that span the tangent space at the singularity for each component, then the kernel is described as~$v_1+v_2+v_3=0$. We computed the action of~$\Ga\rtimes\Gm$ on the tangent space in the description of case \quadric{5}, so it acts by rescaling by~$a_i$, so to preserve the kernel we need~$a_1=a_2=a_3$. The total automorphism group is therefore~$(\Ga^3\rtimes\Gm)\rtimes\Sym_3$.

If we take the Picard group of the conic to be the first term in the N\'eron--Severi group, then the degrees of~$\mathcal{L}_i$ are of the form~$(1,1,0)$ or~$(1,0,1)$. Hence we will reduce~$\Sym_3$ to~$\mathbb{Z}/2\mathbb{Z}$ and then to the trivial group. For the action on~$\Pic^0C$, we replace~$\Ga^3$ by the~$\Ga^2$ given by the kernel of the summation map, and~$\Gm$ acts non-trivially on each non-zero element of~$\Pic^0C$.

\paragraph{Nonreduced, but irreducible}
Only case \quadric{11} is of this form. The automorphism group of the double conic is described in \cite[proposition~8.22]{MR1128218}, and~$\Pic^0C\cong\Ga$. Because~$\PGL_2$ cannot act non-trivially on~$\Pic^0C$ we only need to describe the action of~$\Aut^0C$, the automorphisms of~$C$ which induce the identity on the reduced subscheme, defined by the ideal sheaf~$\mathcal{I}$.

We can interpret the factor~$\Gm$ as automorphisms of the ideal sheaf. Hence these act by rescaling~$\HH^1(C,\mathcal{I})\cong\HH^0(\mathbb{P}^1,\mathcal{O}_{\mathbb{P}^1})$. The factor~$\Ga$ on the other hand coming from~$\Aut^1(C)$ in the notation of~loc.~cit.~will take the degree into account, i.e.~for~$\beta\in\Aut^1(C)$ and~$(n,\lambda)\in\mathbb{Z}\oplus\Pic^0C$ we have
\begin{equation}
  \beta\dot(n,\lambda)=(n,\lambda+n\beta).
\end{equation}
The result now follows.

\paragraph{Nonreduced and reducible}
Cases \quadric{12} and \quadric{13} are of this form. To describe the automorphism group we use the technique from \cite[proposition~8.22]{MR1128218}. For case \quadric{12} the ideal sheaf is the line bundle~$\mathcal{O}_C(-1,-1)$, hence the automorphisms which are the identity on the reduced subscheme are again~$\Ga\rtimes\Gm$. For case \quadric{13} this changes: the ideal sheaf is only supported on the double line, and there it is isomorphic to~$\mathcal{O}(-2)$. The automorphisms which are the identity on the reduced subscheme are therefore only~$\Gm$.

Unlike case \quadric{11} where~$\PGL_2$ acts trivially on~$\Pic^0C$ we have a nontrivial action of the~$\Gm$\dash components coming from automorphisms of the underlying curve. To see this, observe that we can describe~$\Pic^0C$ using the long exact sequence coming from the ``partial normalisation'' which is the morphism from two disjoint double lines to two double lines. We get a long exact sequence
\begin{equation}
  \begin{aligned}
    0&\to k^\times\to(k[x]/(x^2))^\times\oplus(k[y]/(y^2))^\times\to(k[x,y]/(x^2,y^2))^2\to\ldots \\
    &\to\Pic C\to\Pic C_1\oplus\Pic C_2\to 0.
  \end{aligned}
\end{equation}
We realise~$\Pic^0C$ as the cokernel of the appropriate map, which is the vector space spanned by~$xy$. To compute the action on~$\Pic^0C$ by the component of automorphisms coming from the reduced subscheme we again use the description as for case \quadric{5}, and we see that to fix a non-zero element of~$\Pic^0C$ we need the~$\Gm$'s to cancel their action, because the generator~$xy$ gets sent to~$a_1a_2\cdot xy$, where~$a_i\in\Gm$.

The action of~$\Ga\rtimes\Gm$ coming from the automorphisms fixing the reduced subscheme is as in the previous case. Hence we see that in the first step~1~copy of~$\Ga$ is removed, whilst in the second step we remove the copy of~$\Gm$ coming from the automorphisms of the ideal sheaf and we replace~$\Gm^2$ by~$\Gm$.

For case \quadric{12} this means that fixing~$\mathcal{L}_0$ removes the automorphisms which leave the curve and the ideal fixed. In the second step the action of all~$\Gm$'s is nontrivial, and the diagonal subgroup of~$\Gm^2$ is considered. In the third step nothing changes.

For case \quadric{13} on the other hand the degree of the line bundles will play a role. The automorphisms of the ideal sheaf will disappear in the second step as before. Regarding the other components, the factor~$\Sym_2$ disappears in the first step by considering the partial degrees of the line bundle. In the second and third step the~$\Gm$\dash components are paired together in two different ways, so that only a single~$\Gm$ remains.

\begin{table}[p]
  \small
  \renewcommand{\arraystretch}{1.4}
  \centering
  \begin{tabular}{ccccc}
    \toprule
    & decomposition
    & Segre symbol
    & $\Pic^0C$
    & $\Aut C$
    \\
    \midrule
    \quadric{1}
    & $(2,2)$
    & $[1,1,1,1]$
    & $C$
    & $C\rtimes\Aut_p(C)$
    \\
    \quadric{2}
    & $(2,2)$
    & $[3,1]$
    & $\Ga$
    & $\mathbb{Z}/2\mathbb{Z}^{\oplus 2}\rtimes\Aut_p(C)$
    \\
    \quadric{3}
    & $(2,2)$
    & $[2,1,1]$
    & $\Gm$
    & $\Gm\rtimes\Sym_2$
    \\

    \quadric{4}
    & $(1,1)+(1,1)$
    & $[(1,1),1,1]$
    & $\Gm$
    & $\Gm^2\rtimes\Dih_2$ 
    \\
    \quadric{5}
    & $(1,1)+(1,1)$
    & $[(2,1),1]$
    & $\Ga$
    & $(\Ga^2\rtimes\Gm)\rtimes\Sym_2$
    \\
    \quadric{6}
    & $(1,1)+(1,0)+(0,1)$
    & $[2,(1,1)]$
    & $\Gm$
    & $\Gm^3\rtimes\Sym_3$
    \\
    \quadric{7}
    & $(1,1)+(1,0)+(0,1)$
    & $[(3,1)]$
    & $\Ga$
    & $(\Ga^3\rtimes\Gm)\rtimes\Sym_3$
    \\
    \quadric{8}
    & $(1,0)+(1,0)+(0,1)+(0,1)$
    & $[(1,1),(1,1)]$
    & $\Gm$
    & $\Gm^4\rtimes\Dih_4$
    \\
    \quadric{9}
    & $(2,1)+(0,1)$
    & $[2,2]$
    & $\Gm$
    & $\Gm^2\rtimes\Dih_2$
    \\
    \quadric{10}
    & $(2,1)+(0,1)$
    & $[4]$
    & $\Ga$
    & $(\Ga^2\rtimes\Gm)\rtimes\Sym_2$
    \\

    \quadric{11}
    & $2(1,1)$
    & $[(1,1,1),1]$
    & $\Ga$
    & $(\Ga\rtimes\Gm)\times\PGL_2$ 
    \\

    \quadric{12}
    & $2(1,0)+2(0,1)$
    & $[(2,1,1)]$
    & $\Ga$
    & $(\Ga\rtimes\Gm)\times(\Ga\rtimes\Gm)^2\rtimes\Sym_2$
    \\
    \quadric{13}
    & $2(1,0)+(0,1)+(0,1)$
    & $[(2,2)]$
    & $\Ga$
    & $\Gm\times\left( \Gm\times(\Ga\rtimes\Gm)^2 \right)\rtimes\Sym_2$
    \\
    \bottomrule
  \end{tabular}
  \caption{Properties of elliptic quadruples}
  \label{table:hochschild-cubic-properties}
\end{table}

\begin{table}[p]
  \renewcommand{\arraystretch}{1.4}
  \centering
  \begin{tabular}{ccccc}
    \toprule
    & $\Aut_{\mathcal{L}_0}(C)$
    & $\Aut_{\mathcal{L}_0,\mathcal{L}_1}(C)$
    & $\Aut_{\mathcal{L}_0,\mathcal{L}_1,\mathcal{L}_2}(C)$
    \\
    \midrule
    \quadric{1}
    & $\mathbb{Z}/2\mathbb{Z}^{\oplus 2}\rtimes\Aut_p(C)$
    & $\mathbb{Z}/2\mathbb{Z}^{\oplus2}$
    & $\mathbb{Z}/2\mathbb{Z}^{\oplus2}$
    \\
    \quadric{2}
    & $\Gm$
    & 1
    & 1
    \\
    \quadric{3}
    & $\mathrm{V}_4$
    & $\mathbb{Z}/2\mathbb{Z}$
    & $\mathbb{Z}/2\mathbb{Z}$
    \\
    \quadric{4}
    & $\Gm\rtimes\Dih_2$
    & $\Gm\times\mathbb{Z}/2\mathbb{Z}$
    & $\Gm\times\mathbb{Z}/2\mathbb{Z}$
    \\
    \quadric{5}
    & $(\Ga\rtimes\Gm)\times\Sym_2$
    & $\Ga\times\Sym_2$
    & $\Ga\times\Sym_2$
    \\
    \quadric{6}
    & $\Gm^2\rtimes\mathbb{Z}/2\mathbb{Z}$
    & $\Gm$
    & $\Gm$
    \\
    \quadric{7}
    & $(\Ga^2\rtimes\Gm)\rtimes\mathbb{Z}/2\mathbb{Z}$
    & $\Ga^2$
    & $\Ga^2$
    \\
    \quadric{8}
    & $\Gm^3\rtimes\mathbb{Z}/2\mathbb{Z}$
    & $\Gm^2\rtimes\mathbb{Z}/2\mathbb{Z}$
    & $\Gm^2\rtimes\mathbb{Z}/2\mathbb{Z}$
    \\
    \quadric{9}
    & $\Gm\rtimes\mathbb{Z}/2\mathbb{Z}$ or\footnote{Depending on whether you take~$\mathcal{L}_0$ of degree~$(2,0)$ or~$(1,1)$.} $\Gm\rtimes\Dih_2$
    & $\Gm$
    & $\Gm$
    \\
    \quadric{10}
    & $\Ga\rtimes\Gm$ or $(\Ga\rtimes\Gm)\times\Sym_2$
    & $\Ga$
    & $\Ga$
    \\
    \quadric{11}
    & $\Gm\times\PGL_2$
    & $\PGL_2$
    & $\PGL_2$
    \\
    \quadric{12}
    & $\Gm\times(\Ga\rtimes\Gm)^2\rtimes\Sym_2$
    & $\Ga^2\rtimes(\Gm\times\Sym_2)$
    & $\Ga^2\rtimes(\Gm\times\Sym_2)$
    \\
    \quadric{13}
    & $\Gm\times\left( \Gm\times(\Ga\rtimes\Gm)^2 \right)$
    & $(\Ga\rtimes\Gm)^2$
    & $\Ga^2\rtimes\Gm$
    \\
    \bottomrule
  \end{tabular}
  \caption{Automorphism groups of elliptic quadruples}
  \label{table:hochschild-cubic-automorphisms}
\end{table}

\begin{remark}
  \label{remark:translation-principle}
  In \cref{remark:okawa} it was explained that there is a correspondence between equivalences of~$\gr A$ and~$\qgr A$ where~$A$ is a quadratic Artin--Schelter regular~$\mathbb{Z}$\dash algebra. For \emph{cubic} Artin--Schelter regular~$\mathbb{Z}$\dash algebras such a result is \emph{not} true: by the translation principle there exist non-isomorphic~$\mathbb{Z}$\dash algebras whose~$\qgr$'s are equivalent \cite[\S6]{MR2836401}.

  This does not influence the method of proof taken in this paper, but it is interesting to observe that only a discrete amount of information is added, which is irrelevant when taking Lie algebras. It is expected that this is the only new possibility for such a Morita equivalence \cite[\S11]{MR1816070}.
\end{remark}

\bibliographystyle{plain}
\bibliography{mr-clean,arxiv}

\end{document}